\RequirePackage{fix-cm}
\documentclass[smallextended]{svjour3}             
\smartqed  
\usepackage{lineno,hyperref}
\usepackage{graphicx}
\usepackage{csquotes}
\usepackage{epsfig}
\usepackage{lscape}
\usepackage{amsfonts}
\usepackage[misc,geometry]{ifsym}
\usepackage{color}
\usepackage{amsmath,amssymb,dsfont}
\usepackage{listings}
\usepackage{rotating}
\usepackage{multicol,multirow}
\usepackage[english]{babel}
\usepackage{lineno}
\usepackage{enumerate}
\usepackage[justification=centering]{caption}
\usepackage{caption}
\usepackage{subcaption}
\usepackage{afterpage}
\usepackage{empheq}
\modulolinenumbers[5]
\newcommand{\be}{\begin{equation}}
\newcommand{\ee}{\end{equation}}
\newcommand{\bea}{\begin{eqnarray}}
\newcommand{\eea}{\end{eqnarray}}
\newcommand{\nn}{\nonumber}
\newcommand{\bee}{\begin{eqnarray*}}
\newcommand{\eee}{\end{eqnarray*}}

\newcommand{\bt}{\begin{tabbing}}
\newcommand{\et}{\end{tabbing}}
\newcommand{\btb}{\begin{tabular}}
\newcommand{\etb}{\end{tabular}}

\newcommand{\bc}{\begin{center}}
\newcommand{\ec}{\end{center}}

\newtheorem{t1}{Theorem}
\newcommand{\lb}{\label}

\newtheorem{n1}{Note}
\newtheorem{d1}{Definition}
\newtheorem{l1}{Lemma}

\newtheorem{alg1}{Algorithm}

\newcounter{magicrownumbers}
\newcommand\rownumber{\stepcounter{magicrownumbers}\arabic{magicrownumbers}}
\begin{document}

\title{A Trust Region Proximal Gradient Method for Nonlinear Multi-objective Optimization Problems}
\author{Md Abu Talhamainuddin Ansary}    


\institute{Md Abu Talhamainuddin Ansary (\Letter)  \at
              Department of Mathematics,\\ Indian Institute of
              Technology Jodhput,\\ India-342030\\
              \email{md.abutalha2009@gmail.com}}

\date{Received: date / Accepted: date}

\maketitle

\begin{abstract}
In this paper, a globally convergent trust region proximal gradient method is developed for composite multi-objective optimization problems where each objective function can be represented as the sum of a smooth function and a nonsmooth function. The proposed method is free from any kind of priori chosen parameters or ordering information of objective functions. At every iteration of the proposed method, a sub problem is solved to find a suitable direction. This sub problem uses a quadratic approximation of each smooth function and a trust region constraint. An update formula for trust region radius is introduce in this paper. A sequence is generated using descent directions. It is justified that under some mild assumptions every accumulation point of this sequence is a critical point. The proposed method is verified and compared with some existing methods using a set of problems.
\end{abstract}

\section{Introduction}
\label{intro}
In a multi-objective optimization problem, several objective functions are minimized simultaneously. Application of multi-objective optimization can be found in . Classical methods of solving multi-objective optimization problems are scalarization methods (see \cite{kd0,kmm1}), which reduce the original problem to a single objective optimization problem using a set of priori chosen parameters. These methods are user dependent and often fail to generate Pareto front. Heuristic methods like evolutionary algorithms (see \cite{deb2000fast,deb2013evolutionary}), are often used to find approximate Pareto front but can not guarantee any convergence property. Due to these limitations of classical methods, developing gradient/subgradient based methods for multi-objective optimization problems is a major area of research interest. Recently many researchers have developed new techniques for nonlinear multi-objective optimization problems. These methods are possible extension of gradient/subgradient based techniques for single objective optimization to multi-objective case. Gradient based techniques for smooth multi-objective optimization problems are developed in \cite{mat1,mat2,mat4,mat5,nantu1,nantu2,flg1,flg2,flg3,shao2}. Different new techniques are developed in for nonsmooth multi-objective optimization problems in \cite{mat6,bello1,bento1,bonnel1,montonen1,neto1,peng1,tanabe1}. This techniques are possible extension of single objective subgradient methods (\cite{bello1,montonen1,neto1}), proximal point methods (\cite{bento1,bonnel1}), proximal gradient methods (\cite{mat6,peng1,tanabe1}) etc. to multi-objective case.\\ 

Proximal gradient methods are considered as efficient techniques to solve composite single objective optimization problems (see \cite{proxg1,beck2,beck1,nesterov1}). At every iteration of this method linear approximation of smooth function is used to find a suitable descent direction. The ideas of proximal gradient methods are further extended by several researchers in various directions. The proximal gradient method developed by Lee et al. (\cite{lee1}) uses a quadratic approximation of smooth function in every iteration. This method converges quadratically under some mild assumptions.\\

Recently Tanabe et al. (\cite{tanabe1}) have developed a proximal gradient method for multi-objective optimizations. This method combines the ideas of steepest descent method and proximal point method for multi-objective optimization problems developed in \cite{flg2} and \cite{bonnel1} respectively. Similar to the single objective proximal gradient method, the convergence rate of this method is low. The ideas of steepest descent method is replaced by Newton method in \cite{mat6} and by quasi Newton method in \cite{peng1}. In this paper, we have introduced the concept of trust region method and developed a trust region proximal gradient method for composite multi-objective optimization problems.\\

The outline of the paper is as follows. Some prerequisites are discussed in Section \ref{secpre}. A trust region proximal gradient is developed in Section \ref{nprox_method}. An algorithm is proposed in Section \ref{secalg}. The global convergence of this algorithm is justified in this section. In Section \ref{sec_ex}, the proposed method is verified and compared with some existing method using a set of problems.
\section{Preliminaries}
\lb{secpre}
Consider the multi-objective optimization problem:
\bee
(MOP):~~\underset{x\in\mathbb{R}^n}{\min}~~ F(x)=(F_1(x),F_2(x),...,F_m(x)).\\
\eee
Suppose $F_j:\mathbb{R}^n\rightarrow \mathbb{R}$ is defined by $F_j(x)=f_j(x)+g_j(x)$ where $f_j$ is convex and continuously differentiable and $g_j$ is convex and continuous but not necessary differentiable function, for $j=1,2,...,m$. 
Denote $\Lambda_n=\{1,2,...,n\}$ for any $n\in \mathbb{N}$. Inequality in $\mathbb{R}^m$ is understood component wise. If there exists $x\in \mathbb{R}^n$ such that $x$ minimizes all objective functions simultaneously then it is an ideal solution. But in practice, decrease of one objective function may cause increase of another objective function. So in the theory of multi-objective optimization optimality is replaced by efficiency. A point $x^*\in \mathbb{R}^n$ is said to be an efficient solution of the $(MOP)$ if there does not exist $x\in \mathbb{R}^n$ such that $F(x)\leq F(x^*)$ and $F(x)\neq F(x^*)$. A feasible point $x^*\in \mathbb{R}^n$ is said to be a weak efficient solution of the $(MOP)$ if there does not exist $x\in \mathbb{R}^n$ such that $F(x)< F(x^*)$. It is clear that every efficient solution of $(MOP)$ is a weak efficient solution, but the converse is not true. If each $F_j$, $j\in\Lambda_m$ are strictly convex then every weak efficient solution is an efficient solution. For $x, y \in \mathbb{R}^n$, we say $y$ dominates $x$, if and only if $F(y)\leq F(x)$, $ F(y)\neq F(x)$. If $X^*$ is the set of all efficient solutions of the $(MOP)$, then $F(X^*)$ is said to be the Pareto front of the $(MOP)$ and it lies on the boundary of $F(\mathbb{R}^n)$.\\

Let $x^*$ be a weak efficient solution of $(MOP)$. Then $x^*$ must satisfy
\bea
\underset{j\in\Lambda_m}{\max} F_j^{'}(x^*;d)\geq 0 ~~~~\mbox{ for all $d\in\mathbb{R}^n$} \lb{A*}
\eea
The inequality in (\ref{A*}) is sometimes refereed to in the literature as the criticality condition or the first order necessary for weak efficiency of $(MOP)$ and $x^*$ satisfying (\ref{A*}) is often called a critical point for the $(MOP)$. Further convexity of each $F_j$, $j\in \Lambda_m$ ensures that every critical point of $(MOP)$ is weak efficient solution. Further note that if each $F_j$, $j\in\Lambda_m$ is a strictly or strongly convex function then the critical point of $(MOP)$ is an efficient solution. Note that if either $f_j$ or $g_j$ is strictly or strongly convex then is so $F_j$. \\

In nonsmooth optimization, the concept of the gradient (in smooth optimization) is replaced by subdifferential. It plays an important role in nonsmooth optimization. Subdifferential of a convex continuous function is defined as follows.
\begin{d1}\lb{subgrad0}(\cite{jd0})
Suppose $h:\mathbb{R}^n\rightarrow (-\infty,\infty]$ be a proper function and $x\in$ $dom(h)$. Then
subdifferential of $h$ at $x$ is denoted by $\partial h(x)$ and defined as
\bee
\partial h(x):=\{\xi\in\mathbb{R}^n|h(y)\geq h(x)+\xi^T(y-x)\mbox{ for all }y\in \mathbb{R}^n\}.
\eee
If $x\notin dom~h$ then we define $ \partial h(x)=\emptyset$.
\end{d1}
The following properties of subdifferential are often used in the derivation of the proposed methodology.
\begin{t1}\lb{subgrad1}
({\bf Theorem 3.14, \cite{proxg1}}) Let $h:\mathbb{R}^n\rightarrow (-\infty,\infty]$  be a proper convex function, and assume that $x\in int(dom~h)$. Then $ \partial h(x)$ is nonempty and bounded.
\end{t1}
\begin{l1}({\bf Lemma 1, \cite{mat6}}) \lb{lem_mata}
Suppose $0\in Co \underset{ j\in \Lambda_m}{\cup} \partial F_j(x^*)$ for some $x^*\in\mathbb{R}^n$, then $x^*$ is a critical point of the $(MOP)$.
\end{l1}

\section{A trust region proximal gradient method for $(MOP)$}\lb{nprox_method}
In this section a trust region algorithm is developed for $(MOP)$. We construct a sub problem at $x \in\mathbb{R}^n$ to find a direction of $(MOP)$. Using the ideas of trust region method for smooth unconstrained optimization problems, we use the following approximation of $F_j$ at $x.$
\bee
Q_{j}(x,d):=\nabla f_j(x)^Td+\frac{1}{2} d^T B_j(x) d+ g_j(x+d)-g_j(x).
\eee
where $B_j(x)$ is a positive definite approximation of $\nabla^2 f_j(x)$.\\
Define
\bee
Q(x,d):=\underset{j\in\Lambda_m}{\max}~~\nabla f_j(x)^Td+\frac{1}{2} d^T B_j(x) d+ g_j(x+d)-g_j(x).
\eee
Clearly for any fixed $x$, $ Q_{j}$ is continuous $d$ and hence $Q$ is continuous $d$. From Theorem \ref{subgrad1}, $\partial_{d} Q_j(x,d)$ is nonempty and bounded. Denote $$ I(x,d):=\{j\in\Lambda_m|Q(x,d)= Q_{j}(x,d)\}.$$ From Theorems 2.91 and 2.96 of \cite{jd0} we have,
\bee
\partial_{d} Q_j(x,d) &=&\left\{\nabla f_j(x)+B_j(x) d\right\}+\partial_{d} g_j(x+d).\\
\mbox{ and   $ $  }\partial_{d} Q(x,d)&=&Co \underset{j\in I(x,d)}{\bigcup} \partial_{d} Q_j(x,d).\lb{foop1}
\eee
respectively. Finally for any $x\in\mathbb{R}^n$, we solve the following sub problem to find a suitable descent direction of $(MOP)$.
\bee
P(x,\Delta):~~\underset{d\in\mathbb{R}^n}{\min}~~ Q(x,d)& &\\
s.~t.~~~\|d\|^2& \leq & \Delta^2
\eee 
where $\Delta >0$ is the trust region radius. Define $X(\Delta)=\{d\in\mathbb{R}^n|~~ \|d\|^2 \leq \Delta^2\}$.
\begin{n1}
Once can observe that, if $g_j=0$ for all $j\in\Lambda_m$, the $P(x)$ coincides with the sub problem in \cite{shao2}.
\end{n1}
If $B_j(x)$ is positive definite for every $x$ then $Q_j(x,d)$ is strictly convex function in $d$ for every $j\in\Lambda_m$ and $x\in \mathbb{R}^n$. Hence $Q(x,d)$ is strictly convex function in $d$ for every $x\in \mathbb{R}^n.$ This implies, $P(x)$ has a unique finite minimizer.
\\
Denote $d(x,\Delta)=\underset{d\in X(\Delta)}{arg~\min}~ Q(x,d)$ and $t(x,\Delta)=Q(x,d(x)).$ Clearly for every $x\in \mathbb{R}^n,$
\bea
t(x,\Delta)=Q(x,d(x,\Delta))\leq Q(x,0)=0.\lb{tt0} 
\eea
Since $d(x,\Delta)$ is the solution of $P(x,\Delta)$ from Theorem 3.1 of \cite{jd0},
\bee
0\in \partial_{d} Q(x,d(x,\Delta)) + N_{X(\Delta)}(d(x,\Delta))
\eee
From proposition 3.3 in \cite{jd0}, 
\bee
N_{X(\Delta)}(d(x,\Delta))=\left\{\begin{array}{lr}
     \{0\}& if~~d(x,\Delta)\in int(X(\Delta)) \\
     Cone(d(x,\Delta))&otherwise
\end{array}\right.
\eee
\\

Suppose $d(x,\Delta)\in bd(N_{X(\Delta)}(d(x,\Delta)))$. Then from Theorem 2.96 of \cite{jd0}, there exists $\lambda\in\mathbb{R}^{|I(x,d(x,\Delta))|}_{+}$, $\xi_j\in \partial_{d} g_j(x,d(x,\Delta))$ $j\in I(x,d(x,\Delta))$ and $\mu\geq 0$ such that the following conditions hold:
\bee
\underset{j\in I(x,d(x,\Delta))}{\sum}\lambda_j&=&1\\
\underset{j\in I(x,d(x,\Delta))}{\sum} \lambda_j\left(\nabla f_j(x)+B_j(x) d(x,\Delta)+\xi_j \right) +\mu d(x,\Delta)&=&0.
\eee
Substituting $\lambda_j=0$ for all $j\notin I(x,d(x,\Delta))$ and $\mu=0$ for $d(x,\Delta)\in int (X(\Delta))$ we can write
\bea
\underset{j\in \Lambda_m}{\sum}\lambda_j&=&1\lb{kkt1}\\
\underset{j\in \Lambda_m}{\sum} \lambda_j\left(\nabla f_j(x)+B_j  d(x,\Delta)+\xi_j\right)+\mu d(x,\Delta)&=&0\lb{kkt2}\\
\lambda_j\geq 0,~\lambda_j\biggl(\nabla f_j(x)^Td(x,\Delta)+\frac{1}{2} d(x,\Delta)^TB_j(x) d(x,\Delta)\biggr.~~~~~~~\nn\\
\biggl. + g_j(x+d(x,\Delta))-g_j(x)-t(x)\biggr)&=&0~~j\in\Lambda_m\lb{kkt3}\\
\mu \geq 0~~~\mu(\|d(x,\Delta)\|^2-\Delta^2)&=&0 \lb{kkt_tr1} \\ 
\nabla f_j(x)^Td(x,\Delta)+\frac{1}{2} d(x,\Delta)^TB_j(x) d(x,\Delta)+ g_j(x+d(x,\Delta))-g_j(x)&\leq& t(x)~~j\in\Lambda_m\nn\\ \lb{kkt4}. \\
\|d(x,\Delta)\|^2-\Delta^2 &\leq& 0 \lb{kkt_tr2}
\eea
for any $\xi_j\in  \partial_{d} g_j(x+d(x,\Delta))$ for all $j\notin I(x,d(x,\Delta)).$\\

Thus if $d(x,\Delta)$ is the solution of $P(x,\Delta)$ and $t(x,\Delta)=Q(x,d(x,\Delta))$ then there exists $\lambda\in\mathbb{R}^m_{+}$, $\mu\geq 0$ such that $(d(x,\Delta),t(x,\Delta);\lambda,\mu)$ satisfies (\ref{kkt1})-(\ref{kkt_tr2}).
\begin{l1}\lb{lm1}
Suppose $f_j$ is strictly convex function for all $j$. Then $x\in\mathbb{R}^n$ is a critical point of $(MOP)$ if and only if $d(x,\Delta)=0$ for any $\Delta>0$. 
\end{l1}
{\bf Proof:} If possible let $x$ is a critical point of $(MOP)$ and $d(x,\Delta)\neq0$. Since $f_j$ is strictly convex for every $j$, from (\ref{tt0})
 \bea
 \nabla f_j(x)^Td(x,\Delta)+g_j(x+d(x,\Delta))-g_j(x)\leq -\frac{1}{2}d(x,\Delta)^TB_j(x) d(x,\Delta)<0.\lb{lm100}
 \eea
Since $g_j$ is convex, for any $\alpha\in (0,1)$ 
\bea
g_j(x+\alpha d(x,\Delta))-g_j(x)&\leq&\alpha g_j(x+d(x,\Delta))+(1-\alpha)g_j(x)-g_j(x)\nn\\
&=&\alpha\left(g_j(x+d(x,\Delta))-g_j(x)\right). \lb{lm111}
\eea
Using (\ref{lm111}) in (\ref{lm100}),
\bee
& &\alpha \nabla f_j(x)^Td(x,\Delta)+g_j(x+\alpha d(x,\Delta))-g_j(x)\\
&\leq& \alpha \left(\nabla f_j(x)^Td(x,\Delta)+g_j(x+d(x,\Delta))-g_j(x)\right)\\
&<&0.
\eee 
This implies
\bee
\frac{1}{\alpha}\left(\alpha \nabla f_j(x)^Td(x,\Delta)+g_j(x+\alpha d(x,\Delta))-g_j(x)\right)<0
\eee
Taking limit $\alpha\rightarrow 0^{+}$ in the above inequality we have
\bee
\nabla f_j(x)^Td(x,\Delta)+ \underset{\alpha\rightarrow 0^{+}}{\lim} \frac{g_j(x+\alpha d(x,\Delta))-g_j(x)}{\alpha}&<&0\\
i.e. f_j^{'}(x,d(x,\Delta))+g_j^{'}(x,d(x))&<&0
\eee
for all $j\in \Lambda_m$. This implies $\underset{j\in\Lambda_m}{\max}~F_{j}^{\prime}(x,d(x,\Delta))<0$. This shows that $x$ is not a critical point, a contradiction. Hence if $x$ is a critical point then $d(x,\Delta)=0.$\\

Conversely suppose $d(x,\Delta)=0$ for any $\Delta>0$.\\Then $N_{X(\Delta)} (d(x,\Delta))=\{0\}$ holds since $d(x,\Delta)\in int X(\Delta)$. Hence $\mu=0$ holds in (\ref{kkt_tr1}) and (\ref{kkt_tr2}). Then  from (\ref{kkt1}) and (\ref{kkt_tr2}), there exists $\lambda\in\mathbb{R}^m_{+}$ such that $\sum_{j\in\Lambda_m} \lambda_j=1$ and $\sum_{j\in \Lambda_m} \lambda_j\left(\nabla f_j(x)+\xi_j\right)=0$ where $\xi_j\in\partial g_j(x)$ for $j\in\Lambda_m$. This implies $$0\in Co \underset{ j\in \Lambda_m}{\cup} \partial F_j(x).$$       
Hence from  Lemma \ref{lem_mata}, $x$ is a critical point of $(MOP).$\qed
\begin{n1}
Since $P(x,\Delta)$ has unique solution, $t(x,\Delta)=0$ holds if and only if $d(x,\Delta)=0$. Hence from Lemma \ref{lm1} and (\ref{tt0}) we can conclude that $t(x,\Delta)<0$ holds if and only if $x$ is a non-critical point of $(MOP)$.
\end{n1}
\begin{t1}\lb{t1}
Suppose $d^TB_j(x)d\geq \sigma \|d\|^2$ holds for every $x,d\in\mathbb{R}^n$ and $j\in\Lambda_m$. Further suppose, $d(x,\Delta)$ is the optimal solution of $P(x,\Delta)$. Then 
\bea
t(x,\Delta)\leq-(\frac{\sigma}{2}+\mu) \|d(x,\Delta)\|^2\lb{t10}
\eea 
and $d(x,\Delta)$ is a descent direction for every $F_j$.
\end{t1}
\textbf{Proof:} Suppose $d(x,\Delta)$ is the solution of $P(x,\Delta)$ and $t(x,\Delta)=Q(x,d(x,\Delta))$. Then there exists $\lambda\in\mathbb{R}^m_{+}$ and $\mu\geq 0$ such that $(d(x,\Delta),t(x,\Delta);\lambda,\mu)$ satisfies (\ref{kkt1})-(\ref{kkt_tr2}). Since $g_j$ is convex and $\xi_j\in\partial_{d}g_j(x+d(x,\Delta)),$
\bea
g_j(x+d(x,\Delta))-g_j(x)\leq \xi_j^T d(x,\Delta).\lb{t11}
\eea
Multiplying both sides of (\ref{kkt2}) by $d(x,\Delta)$,
\bee
\underset{j\in \Lambda_m}{\sum} \lambda_j\left[\nabla f_j(x)^Td(x,\Delta)+d(x,\Delta)^TB_j(x) d(x,\Delta)+\xi_j^Td(x,\Delta)\right]+\mu d(x,\Delta)^Td(x,\Delta)=0
\eee
Hence from (\ref{t11}),
\bea
\underset{j\in \Lambda_m}{\sum} \lambda_j\left\{\nabla f_j(x)^Td(x,\Delta)+d(x,\Delta)^TB_j(x) d(x,\Delta)+g_j(x+d(x,\Delta))-g_j(x)\right\}+\mu d(x,\Delta)^Td(x,\Delta)\leq 0.\nn\\ \lb{t12}
\eea
Taking sum over $j\in\Lambda_m$ in (\ref{kkt3})and using (\ref{kkt1}),
\bea
& &\underset{j\in\Lambda_m}{\sum}\lambda_j\left\{\nabla f_j(x)^Td(x,\Delta)+ d(x,\Delta)^TB_j(x) d(x,\Delta)+ g_j(x+d(x,\Delta))-g_j(x)\right\}\nn\\
&= &\underset{j\in\Lambda_m}{\sum}\lambda_j\frac{1}{2}d(x,\Delta)^TB_j(x) d(x,\Delta)+t(x,\Delta).\lb{t13}
\eea 
Using (\ref{t12}) in (\ref{t13}),
\bea
t(x,\Delta)\leq - \underset{j\in\Lambda_m}{\sum}\lambda_j\frac{1}{2}d(x,\Delta)^TB_j(x) d(x,\Delta)-\mu d(x,\Delta)^Td(x,\Delta).\lb{t14}
\eea 
Since $d(x,\Delta)^TB_j(x) d(x,\Delta)\geq \sigma \|d(x,\Delta)\|^2$ holds for every $j$, from (\ref{t14}) and (\ref{kkt1}),
\bee
t(x,\Delta)\leq-(\frac{\sigma}{2}+\mu) \|d(x,\Delta)\|^2. 
\eee
Above inequality shows if $d(x,\Delta)\neq 0$ then $t(,\Delta)<0$. Then from (\ref{kkt4}), 
\bee
\nabla f_j(x)^Td(x,\Delta)+g_j(x+d(x,\Delta))-g_j(x)&\leq& t(x,\Delta)-\frac{1}{2}d(x,\Delta)^TB_j(x)d(x,\Delta)\\
&\leq& -(\sigma+\mu)\|d(x,\Delta)\|^2\\&<&0
\eee
holds for every $j\in\Lambda_m$. This implies $d(x,\Delta)$ is a descent direction for every $f_j$.\qed
\subsection{Algorithm}\lb{secalg}
In this section we develop an algorithm for $(MOP)$ using the theoretical results developed so far. With some initial approximation $x^0$ and initial trust region radius $\Delta_0$ a sequence is generated by update formula $x^{k+1}=x^k+d(x^k,\Delta_k)$, where $d(x^k,\Delta_k)$ is the optimal solution of $P(x^k,\Delta_k)$. For simplicity, rest of the paper $d(x^k,\Delta_k)$ and $t(x^k,\Delta_k)$ are denoted by $d^k$ and $t^k$ respectively.\\

Since we are not using any line search technique, proper selection of trust region radius ($\Delta_k$) is necessary in this algorithm. We follow the following strategies to update $\Delta_{k+1}$. Suppose actual reduction using $d^k$ is $$Ared(d^k)=\underset{j\in\Lambda_m}{\min} \{F_j(x^k)-F_j(x^k+d^k)\}.$$ Predicted reduction using $d^k$ is defined by
 $$Pred(d^k)=-\underset{j\in\Lambda_m}{\max} \{{\nabla f_j(x^k)}^T d^k+\frac{1}{2} {d^k}^TB_j(x^k)d^k+g_j(x^k+d^k)-g_j(x^k)\}=-t^k.$$
 Let $\rho(d^k)=\frac{Ared(d^k)}{Pred(d^k)}$. Under following conditions $\Delta_{k+1}$ is either expanded or shrunk or kept unchanged in next iteration with the help of scalars\\$0<\sigma_3<1<\sigma_1$, $0<\sigma_0< \sigma_2<1$.
 \begin{itemize}
\item $\rho(d^k)\geq \sigma_2$ implies good agreement between $Ared (d^k)$ and $Pred (d^k)$. In this case, trust region radius can be expended for next iteration. Here trust region radius is updated by $\Delta_{k+1}=\max \{\sigma_1\Delta_k,\Delta_{min}\}$ where $\Delta_{min}$ is pre-specified.
\item $\sigma_0\leq \rho(d^k)<\sigma_2$ implies $d^k$ is descent for every $j$ but agreement between $Ared (d^k)$ and $Pred (d^k)$ is not good. In this case, trust region radius remain unchanged.
\item $ \rho(d^k)<\sigma_0$ implies, $d^k$ does not provide sufficient decrease in at least one $j$. In this case, trust region radius is shrunked by $\Delta_k=\sigma_3 \Delta_k  $ trust region and solved $P(x^k,\Delta_k)$ again. This process is repeated until $ \rho(d^k)\geq\sigma_0$ holds.
\end{itemize}
We generate next iterating point $x^{k+1}=x^k+d^k$. For $x^{k+1}$, $B_j(x^{k+1})$ is updated using modified BFGS updated formula in \cite{powell}. This process is repeated until we get an approximate critical point $(MOP)$. Above steps are presented in the following algorithm.
\begin{alg1} \lb{trpxm1}(Trust region proximal gradient method for $(MOP)$)
\begin{enumerate}[{Step} 1]
\item Choose initial approximation $x^0$, initial positive definite matrix $B_j(x^0)$ for $j\in\Lambda_m$, scalars $0<\sigma_3<1<\sigma_1$, $0<\sigma_0< \sigma_2<1$, initial trust region radius $\Delta_0$ and $\Delta_{min}>0$.   Set $k:=0$
\item Solve the sub problem $P(x^k,\Delta_k)$ to find $d^k$. Then compute $t^k=Q(x^k,d^k)$.\lb{P_S}
\item If $\|d^k\|<\epsilon$, then stop. Else go to Step \ref{St_inn}. \lb{st_term}
\item \lb{St_inn} Compute 
\bea
\rho(d^k)=\frac{\underset{j\in\Lambda_m}{\min}~\{F_j(x^k)-F_j(x^k+d^k)\}}{-t^k} \lb{tr_up}
\eea
If $\rho(d^k)<\sigma_0$ set 
\bea
\Delta_k=\sigma_3 \Delta_k  \lb{tr_up1}
\eea
and go to Step \ref{P_S}. If $\rho(d^k)\geq \sigma_0$ then set
\bea
x^{k+1}&=&x^k+d^k  \lb{tr_up2} \\
\Delta_{k+1}&=&\left\{\begin{array}{lr}
\max \{\sigma_1\Delta_k,\Delta_{min}\} & \text{ if } \rho(d^k)\geq \sigma_2 \\
\Delta_k & \text{otherwise}
\end{array} \right.
\eea
\item Generate $B_j(x^{k+1})$ for every $j\in\Lambda_m$ using modified BFGS update formula in \cite{powell}. set $k:=k+1$, and go to Step \ref{P_S}.
\end{enumerate}
\end{alg1}
\section{Convergence analysis of Algorithm \ref{trpxm1}}\lb{sec_conv}
In this section, global convergence of this algorithm is justified under certain assumptions. Prior to that we discuss some properties of $\theta(x; \Delta)$ defined as
\bea
\theta(x; \Delta)=\underset{\|d\|\leq \Delta}{\min}~\left[\underset{j\in\Lambda_m}{\max}~\left\{\nabla f_j(x)^Td+g_j(x+d)-g_j(x)\right\}\right].
\lb{thxd}
\eea
\begin{l1}
Following results hold for $\theta(x; \Delta)$ defined as (\ref{thxd}). Proof of these results are straightforward hence not explained here. 
\begin{enumerate}[(i)]
    \item $\theta(x; \Delta)\leq 0$ for every $x\in\mathbb{R}^n$ and $\Delta>0$.
    \item $\theta(x; \bar{\Delta})\leq \theta(x;\Delta)$  for any $\bar{\Delta} \geq \Delta$.
    \item $x^*$ is a critical point of $(MOP)$ if and only if $\theta(x;\Delta)=0$ for some $\Delta>0$.
    \item $\theta(x;\alpha \Delta)\leq \alpha \theta(x;\Delta)$ for any $\alpha\in[0,1]$.
\end{enumerate}
\end{l1}
\begin{l1} \lb{lm_m1}
Suppose there exist $b>0$ such that $d^TB_{j}(x^k) d \leq b\|d\|^2$ holds for every $j\in\Lambda_m$ and $d\in\mathbb{R}^n.$ Then for any $\Delta\geq \Delta_k$,
\bea
t^k \leq \frac{\theta(x^k;\Delta)}{2\Delta} \min \left\{\Delta_k, \frac{-\theta(x^k;\Delta)}{b\Delta}\right\}\lb{dt1}
\eea
\end{l1}
{\bf Proof:} Suppose $\bar{d}$ is the minima of $\theta(x^k,\Delta_k)$ in (\ref{thxd}). Then
\bea
t(x^k;\Delta_k)&=& \underset{\|d\|\leq \Delta_k}{\min}~\underset{j\in\Lambda_m}{\max}~~\left\{\nabla f_j(x^k)^Td+\frac{1}{2} d^TB_j(x^k)d+g_j(x^k+d)-g_j(x^k)\right\} \nn\\
& \leq & \underset{0\leq \alpha\leq 1}{\min}~\underset{j\in\Lambda_m}{\max}~~\left\{\nabla f_j(x^k)^T(\alpha\bar{d})+\frac{\alpha^2}{2} {\bar{d}}^TB_j(x^k){\bar{d}}+g_j(x^k+\alpha{\bar{d}})-g_j(x^k)\right\} \nn \\
&\leq& \underset{0\leq \alpha\leq 1}{\min}~\underset{j\in\Lambda_m}{\max}~~\left\{\alpha\left[\nabla f_j(x^k)^T\bar{d}+g_j(x^k+\bar{d})-g_j(x^k)\right]+\frac{b\alpha^2}{2} \|\bar{d}\|^2 \right\} \nn\\
&\leq &\underset{0\leq \alpha\leq 1}{\min} \left\{ \alpha \theta(x^k;\Delta_k)+\frac{b\alpha^2}{2} \|\bar{d}\|^2 \right\} \nn\\
&\leq &  \max\left\{\frac{\theta(x^k;\Delta_k)}{2}, \frac{-\theta(x^k;\Delta_k)^2}{2b\|\bar{d}\|^2}\right\} \lb{thxd1}
\eea
Now for any $\Delta>\Delta_k$,
\bea
\theta(x^k;\Delta_k)&=&\theta(x^k;\Delta\frac{\Delta_k}{\Delta})\leq \frac{\Delta_k}{\Delta}\theta(x^k;\Delta)  \lb{thxd2a}
\eea
Since $\theta(x^k;\Delta_k)$ and $\theta(x^k;\Delta)$ are negative, from (\ref{thxd2a}) we have $$\theta(x^k;\Delta_k)^2\geq  \frac{\Delta_k^2}{\Delta^2}\theta(x^k;\Delta)^2.$$ This along with $\|\bar{d}\|\leq \Delta_k$ imply
\bea
\frac{-\theta(x^k;\Delta_k)^2}{2b\|\bar{d}\|^2}\leq \frac{-\theta(x^k;\Delta)^2}{2b\Delta^2} \lb{thxd2b}
\eea
Using (\ref{thxd2a}) and (\ref{thxd2b}) in (\ref{thxd1}),
\bee
t(x^k;\Delta_k)&\leq& \max\left\{ \frac{\Delta_k}{2\Delta}\theta(x^k;\Delta) ,  \frac{-\theta(x^k;\Delta)^2}{2b\Delta^2}\right\}\\
&=&\frac{\theta(x^k;\Delta)}{2\Delta} \min\left\{\Delta_k, \frac{-\theta(x^k;\Delta)}{b\Delta}\right\}
\eee
Last equality holds since $\theta(x^k;\Delta)<0$. Hence the lemma follows.\qed

In the following lemma, we show that inner circle (Step \ref{P_S}- Step \ref{St_inn} -Step \ref{P_S}) is repeated finitely.
\begin{l1}
Suppose assumptions of Lemma \ref{lm_m1} are true. Then for any non critical point $x^k$ $\rho_k\geq \sigma_0$ holds after finite repetition of inner circle (Step \ref{P_S}- Step \ref{St_inn} -Step \ref{P_S}) in Algorithm \ref{trpxm1}.
\end{l1}
{\bf Proof:} If possible suppose the inner circle inner circle (Step \ref{P_S}- Step \ref{St_inn} -Step \ref{P_S}) is repeated infinitely a non critical point $x^k$. Then
\bea
\rho_k(d^{k,r})< \sigma_0~\mbox{and}~\Delta_{k,r+1}=\sigma_3 \Delta_{k,r} \lb{incr1}
\eea
hold for $r=1,2,\dots$, where index $j$ denotes the number of repetition of inner circle and $d^{k,r}$ is the solution of $SP(x^k,\Delta_{k,r})$ .
Since $0<\sigma_3<1$, (\ref{incr1}) implies
\bea
\underset{r\rightarrow \infty}{\lim}  \Delta_{k,r} =0 \mbox{and}~\underset{r\rightarrow \infty}{\lim} \|d^{k,r}\|=0 \lb{incr2}
\eea
This together with (\ref{dt1}) imply that there exist $\bar{j}$ and $\Delta>0$ such that for all $j\geq \bar{j}$ and $\Delta\geq \Delta_{k,j} $,
\bea
t^{k,r}\leq \frac{\theta(x^k,\Delta))\Delta_{k,r}}{\Delta} \leq -\delta \Delta_{k,r} \lb{incr3}
\eea
where $t^{k,r}:=t(x^k,\Delta_{k,r})$.
Since $B_j(x^{k})\approx\nabla^2f_j(x^k)$, there exists $\omega_j(x^k;d)$ such that
\bee
f_j(x^k+d^{k,r})=f_j(x^k)+\nabla f_j(x^k)^Td^{k,r}+\frac{1}{2} {d^{k,r}}^TB_j(x^k)d^{k,r}+\omega_j(x^k;d^{k,r})
\eee
and $\frac{\omega_j(x^k;d)}{\|d\|^2}\longrightarrow 0$ as $\|d\|^2\longrightarrow 0$ for all $j$. Above equality implies 
 \bee
f_j(x^k+d^{k,r})+g_j(x^k+d^{k,r})-f_j(x^k)-g_j(x^k)&=&\nabla f_j(x^k)^Td^{k,r}+\frac{1}{2} {d^{k,r}}^TB_j(x^{k,r})d^{k,r}\\& &+g_j(x^k+d^{k,r})-g_j(x^k)+\omega_j(x^k;d^{k,r})
\eee
Define $\omega (x^k;d)=\underset{j\in\Lambda_m}{\max}~\omega_j(x^k;d)$. Clearly $\frac{\omega(x^k;d)}{\|d\|^2}\longrightarrow 0$ as $\|d\|^2\longrightarrow 0$ and 
\bee
\underset{j\in\Lambda_m}{\max} F_j(x^k+d^{k,r})-F_j(x^k) &\leq & \underset{j\in\Lambda_m}{\max} \left\{\nabla f_j(x^k)^Td^{k,r}+\frac{1}{2} {d^{k,r}}^TB_j(x^k)d^{k,r}+g_j(x^k+d^{k,r})-g_j(x^k) \right\}+\omega(x^k;d^{k,r})
\eee
This implies 
\bee
-\underset{j\in\Lambda_m}{\min} F_j(x^k)-F_j(x^k+d^{k,r}) -t(x^{k,r}) \leq  \omega(x^k;d^{k,r})
\eee
Now $\rho(d^{k,r})<1$ implies $ \underset{j\in\Lambda_m}{min}~~\left\{ F_j(x^k)-F_j(x^k+d^{k,r})\right\} +t(x^{k,r})<0$. Hence from above inequality
\bee
|\underset{j\in\Lambda_m}{\min} F_j(x^k)-F_j(x^k+d^{k,r})+t(x^{k,r})| \leq \omega(x^k;d^{k,r})  
\eee
This along with (\ref{incr3}) implies 
\bee
\frac{\left| \underset{j\in\Lambda_m}{min}~~\left\{ F_j(x^k)-F_j(x^k+d^{k,r})\right\} +t(x^{k,r})\right|}{-t(x^{k,r})}\leq \frac{\omega(x^k;d^{k,r})}{\delta \Delta_{k,r}}
\eee
Since $t(x^{k,r})<0$ and $\|d^{k,r}\|\leq \Delta_{k,r}$, above inequality implies
\bea
\left| \rho(d^{k,r}) -1 \right|\leq \frac{1}{\delta} \frac{\omega(x^k;d^{k,r})}{\|d^{k,r}\|^2} \frac{\|d^{k,r}\|^2}{\Delta_{k,r}^2} \Delta_{k,r}\leq \frac{1}{\delta} \frac{\omega(x^k;d^{k,r})}{\|d^{k,r}\|^2} \Delta_{k,r} \lb{dttr}
\eea
Since $\underset{r \rightarrow \infty}{\lim}~\frac{\omega(x^k;d^{k,r})}{\|d^{k,r}\|^2} =0$ and $\underset{r \rightarrow \infty}{\lim}\Delta_{k,r}=0$ there exist $0<\delta_1<1$  such that $ \frac{1}{\delta} \frac{\omega(x^k;d^{k,r})}{\|d^{k,r}\|^2}  \Delta_{k,r} \leq \delta_1$ holds for every $r$ sufficiently large. Then (\ref{dttr}) implies, $\rho(d^{k,r})\geq 1-\delta_1$ holds for every $r$ sufficiently large. This contradicts that inner circle is repeated infinitely.\qed
$ $\\

Finally we justify the global convergence of Algorithm \ref{trpxm1} in the following theorem under some mild assumptions.
\begin{t1}\lb{cv0}
Suppose $\{x^k\}$ be a sequence generated by Algorithm \ref{trpxm1}. Further suppose assumptions of Lemma \ref{lm_m1} are true and the level set\\ $L=\{x:F(x)\leq F(x^0)\}$ is bounded. Then either $\{x^k\}$ terminates at a critical point of $(MOP)$ or any accumulation point of $\{x^k\}$ is a critical point of $(MOP)$. 
\end{t1}
{\bf Proof:} Since inner circle (Step \ref{P_S}- Step \ref{St_inn} -Step \ref{P_S}) is repeated finitely, Algorithm \ref{trpxm1} terminates at $x^k$ if and only if $d^k=0$. Then from (\ref{kkt_tr1}), $\mu_k=0$. Substituting $d^k=0$ and $\mu_k=0$ in (\ref{kkt1}) and (\ref{kkt2}), we have
\bee
\underset{j\in \Lambda_m}{\sum} \lambda_j^k\left(\nabla f_j^k+\xi_j^k\right)=0
\eee
where ${\xi}_j^k\in \partial g_j(x^k)$, $\lambda_j^k\geq 0$, and $\underset{j\in \Lambda_m}{\sum}\lambda_j=1$. This implies  $0\in Co\left(\underset{j\in\Lambda_m}{\bigcup} \partial F_j(x^k)\right)$. Hence from Lemma \ref{lem_mata}, $x^k$ is a critical point of $(MOP)$.\\

Next suppose Algorithm \ref{trpxm1} generates an infinite sequence. Since $r(d^k)\geq \sigma_0$ holds for every $k$, we have
\bea
-\sigma_0 t^k&\leq& \underset{j\in\Lambda_m}{\min} \left\{F_j^k-F_j^{k+1}\right\} \nn\\
&\leq&  \underset{j\in\Lambda_m}{\min} \left\{F_j^{k-1}-F_j^{k+1}\right\}- \underset{j\in\Lambda_m}{\min} \left\{F_j^{k-1}-F_j^{k}\right\} \lb{min_ineq} \\
&\leq& \underset{j\in\Lambda_m}{\min} \left\{F_j^{k-1}-F_j^{k+1}\right\} +\sigma_0t^{k-1} \nn
\eea
This implies $-\sigma_0 t^k-\sigma_0t^{k-1}\leq \underset{j\in\Lambda_m}{\min} \left\{F_j^{k-1}-F_j^{k+1}\right\} $. Proceeding by similar way, we have
\bee
\sigma_0\sum_{i=0}^k (-t^i) \leq  \underset{j\in\Lambda_m}{\min} \left\{F_j^{0}-F_j^{k+1}\right\}
\eee
Since $x^k\in L$ for all $k$, $x^k$ is bounded. Then there exists $F_j^*>-\infty$ for all $j$ such that  
\bee
\sigma_0\sum_{i=0}^{\infty} (-t^i) \leq  \underset{j\in\Lambda_m}{\min} \left\{F_j^{0}-F_j^*\right\} <\infty
\eee
This implies $t^k\rightarrow 0 $ as $k\rightarrow \infty$ as $-t^k>0$ for all $k$. Then from (\ref{t10}),\\ $\|d^k\|\rightarrow 0$ as $k\rightarrow \infty$, since $\mu^k\geq 0$ for all $k$.\\

Since $\{x^k\}$ is bounded, there exists at least one convergent sub sequence. Suppose  $\{x^k\}_{k\in K_0}$ be a converging subsequence of $x^k$, converging to $x^*$. Next we show that, $x^*$ is a critical point of $(MOP)$.\\

 Since $(t(x^k,\Delta_k),d(x^k,\Delta_k))$ be the solution of $P(x^k,\Delta_k)$, there exists $(\lambda^k,\mu_k)$ satisfying (\ref{kkt1})-(\ref{kkt_tr2}). From (\ref{kkt1}), $\{\lambda^k\}$ is bounded. So there exists a subsequence of $\{\lambda^k\}$ converging to $\lambda^*$. Without loss of generality, we may assume $\lambda^k\rightarrow \lambda^*$ as $k\rightarrow \infty$. Then $\lambda_j^*\geq 0$ for all $j$ and $\sum_{j\in\Lambda_m} \lambda_j^{*}=1$ (from (\ref{kkt1})). Now $d(x^k)\rightarrow 0$ and $\Delta_k\geq \Delta>0$ imply, $\mu_k=0$ for large $k$. Similarly $\|B_j(x^k)\| \leq b$ implies $\|B_j(x^k)d^k\|\leq \|B_j(x^k)\|\|d^k\|\leq b\|d^k\|$. This along with $d^k\rightarrow \infty$ as $k\rightarrow \infty$ implies $B_j(x^k)d^k\rightarrow 0$ as $k\rightarrow \infty$. Thus taking limit $k\rightarrow \infty$ in (\ref{kkt1}) and (\ref{kkt2}), we have 
\bee
\underset{j\in\Lambda_m}{\sum} \lambda_j^*&=&1\\
\underset{j\in\Lambda_m}{\sum} \lambda_j^*(\nabla f_j(x^*)+\xi_j^*)&=&0
\eee
This implies $0\in Co\left(\underset{j\in\Lambda_m}{\bigcup} \partial F_j(x^*)\right).$ Hence from Lemma \ref{lem_mata}, $x^*$ is a critical point of $(MOP)$. Since $\{x^k\}_{k\in K_0}$ is an arbitrary convergent subsequence of $\{x^k\}$,  any accumulation point of $\{x^k\}$ is a critical point of $(MOP)$.\qed
\begin{n1}
In some articles in literature (\cite{nantu1,shao2}), actual reduction using $d^k$ is defined as $Ared(d^k)=\underset{j\in\Lambda_m}{\max} \{f_j(x^k)-f_j(x^k+d^k)\}.$ In this case $\rho(d^k)\geq \sigma_2$ does not necessarily imply sufficient decrease in each objective function. In addition to this, Inequality (\ref{min_ineq}) does not holds for maximum. That is why, we have defined actual reduction different from \cite{nantu1,shao2}.
\end{n1}
\section{Numerical examples}
\lb{sec_ex}
In this section Algorithm \ref{trpxm1} (MOTRPG) is verified and compared with proximal gradient method for $(MOP)$ developed in \cite{mat6} (MONPG) and \cite{tanabe1} (MOPG) using a set of problems. Implementation of these algorithms are explained below.
\begin{itemize}
\item MATLAB (2024a) code is developed for each method using the extension `{\it CVX}'. `{\it CVX}' programming with solver `{\it SeDuMi}' is used to solve the sub problems in each method.
\item In MOTRPG, we have used $\sigma_0=0.01$, $\sigma_1=1.5$, $\sigma_2=0.5$, $\sigma_3=0.5$, and $\Delta_{min}=\max\{\underset{j\in\Lambda_m}{\min} \|f_j(x^0)\|, 1\}.$
\item $\|d^k\|<10^{-5}$ or maximum 2000 iterations is considered as stopping criteria.
\item Solution of a multi-objective optimization problem is not isolated optimum points but a set of efficient solutions. To generate an approximate set of efficient solutions we have considered multi-start technique. Following steps are executed for each method.
\begin{itemize}
\item A set of 100 uniformly distributed random initial points between $lb$ and $ub$ are considered, where $lb,ub\in\mathbb{R}^n$ and $lb<ub$.
\item Algorithm of each method is executed individually.
\item Suppose $\mathcal{WX^*}$ is the collection of approximate critical points. The non dominated set of $\mathcal{WX^*}$ is considered as an approximate set of efficient solutions.
\end{itemize} 
\end{itemize}
Next we explain the steps of Algorithm \ref{trpxm1} using the following example. \\
\textbf{Example 1:} Consider the bi-objective optimization problem: 
$$(E_1):~~\underset{x\in\mathbb{R}^2}{\min}~~\left(f_1(x)+g_1(x), f_2(x)+g_2(x)\right)$$
where
\bea
f_1(x)&:=& x_1^2+x_2^2,\\
f_2(x)&:=&(x_1-5)^2+(x_2-5)^2,\\
g_1(x)&:=& \max\left\{(x_1-2)^2+(x_2+2)^2,x_1^2+8x_2\right\},\\
g_2(x)&:=& \max\left\{5x_1+x_2,x_1^2+x_2^2\right\}.
\eea
Here $f_1,f_2:\mathbb{R}^2\rightarrow \mathbb{R}$ follows from test problem BK1 in and $f_1,f_2:\mathbb{R}^2\rightarrow \mathbb{R}$ follows from Example 2 in Section 4.1 of \cite{montonen1}.\\\\
Consider $x^0=(-4.5,6.5 )^T.$ Then $F(x^0)=(177,155)^T$.\\ Using $\Delta_0=\min\{\|\nabla f_1(x^0)\|^2,\|\nabla f_1(x^0)\|^2\}=15.8114$, solution of $P(x^0;\Delta_0)$ is obtained as $d^0=(3.4524,-1.9728)^T$ and $t^0=-104.5165$.\\ Then $\rho(d^0)= 0.9244>\sigma_2$. So we calculate next iterating point\\ $x^1=x^0+d^0=( 1.1667,-1.8333)^T$ and expand the trust region for next iteration as $\Delta_1=\max\{\sigma_3\Delta_0,\Delta_{min}\}=18.9737.$ One can observe that\\$F(x^1)=(73.4843,58.3892)^T<F(x^0)$. Using the stopping criteria\\ $\|d^k\|<10^{-4}$, approximate solution is obtained as$$x^3=( 2.0077,3.0015)^T\approx(2,3)^T.$$

\begin{figure}[!htbp]
   \centering
      \includegraphics[height=2.75cm,width=.95\linewidth]{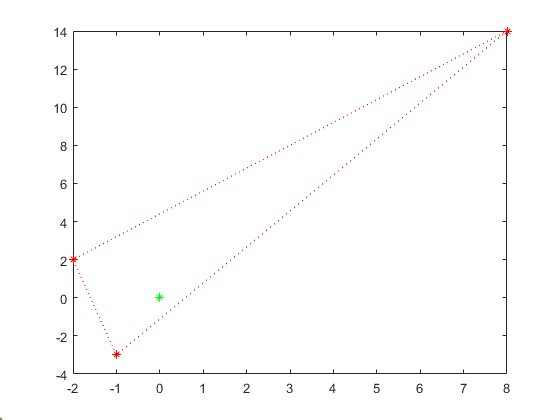}
       \caption{Convex hull of $\left\{\partial F_1(x^*),\partial F_2(x^*)\right\}$}
       \label{fig0}
\end{figure}

Now we justify that $x^*=(2,3)^T$  is a critical point of $(MOP)$. One can observe that both $x_1^2+8x_2$ is active function in $g_1(x^*)$ and both $5x_1+x_2$ and $x_1^2+x_2^2$ are active functions in $g_2(x^*)$. Hence from Theorem 2.96 of \cite{jd0},\\ $\partial g_1(x^*)=\{(4,8)^T\}$ and $\partial g_2(x^*)=Co\{(5,1)^T, (4,6)^T\}.$ This implies 
\bee
 \partial F_1(x^*)&=&\{\nabla f_1(x^*)\}+\partial g_1(x^*)\\&=& \{(4,6)^T\}+\{(4,8)^T\}\\
 &=&\{(8,14)T\} \\
 \partial F_2(x^*)&=&\{\nabla f_2(x^*)\}+\partial g_2(x^*)\\&=&\{(-6,-4)^T\}+C0\{(5,1)^T, (4,6)^T\}\\&=&Co\{(-1, -3)^T, (-2,2)^T\}
 \eee
 Then convex hull of $\{\partial F_1(x^*),\partial F_2(x^*)\}$ is the triangle formed by $(8,14)$, $(-1,-3)$, and $(-2,2)$. From Figure \ref{fig0}, we can see that $(0,0)$ lies inside the triangle.
 $ $\\
 Hence we can say $0\in Co\{\partial F_1(x^*),\partial F_2(x^*)\}$. This implies $x*=(2,3)^T$ is a critical point of Example 1. \\

\begin{figure}[!htbp]
    \centering
       \includegraphics[height=2.75cm,width=.95\linewidth]{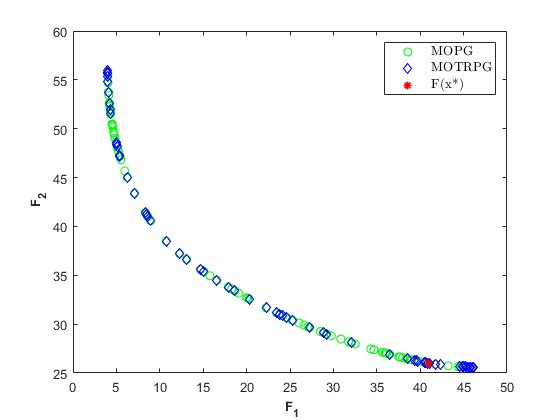}
        \caption{Approximate Pareto fronts of $P_1$}
        \label{fig1}
\end{figure}

 Approximate Pareto front of $(E_1)$ obtained using multi-start technique with MONPG and MOPG is provided in Figure \ref{fig1}. One can observe that $F(x^*)$ belongs to this approximate Pareto front.
$ $ \\

{\bf $l_1$-regularized multi-objective least square problem (MOLS):} Least square problems play an important role different fields of data science. A multi-objective least square problem is of the form
$$\underset{x\in\mathbb{R}^n}{\min}~~\left( \frac{1}{2}\|A_1x-b_1\|^2, \frac{1}{2}\|A_2x-b_2\|^2, \dots,\frac{1}{2}\|A_mx-b_m\|^2\right)  $$
where $A_j\in \mathbb{R}^{m\times n}$ and $b_j\in  \mathbb{R}^{m}$ for all $j\in\Lambda_m$.
However, many least square problems seek a sparse solution. In such cases a regularized version on the least square problem is solved by adding the
a multiple of the $l_1$ norm of $x$ to the main objective. The regularized multi-objective least square problem can be given as
$$\underset{x\in\mathbb{R}^n}{\min}~~\left( \frac{1}{2}\|A_1x-b_1\|^2+ \frac{\nu_1}{2} \|x\|_1 , \frac{1}{2}\|A_2x-b_2\|^2+\frac{\nu_2}{2} \|x\|_1, \dots,\frac{1}{2}\|A_mx-b_m\|^2+\frac{\nu_m}{2} \|x\|_1\right)  $$
where $\mu_j>0$ are regularized parameters. For details, readers may see \cite{boyd2}\\

\begin{figure}[!htbp]
    \centering
       \includegraphics[height=2.75cm,width=.95\linewidth]{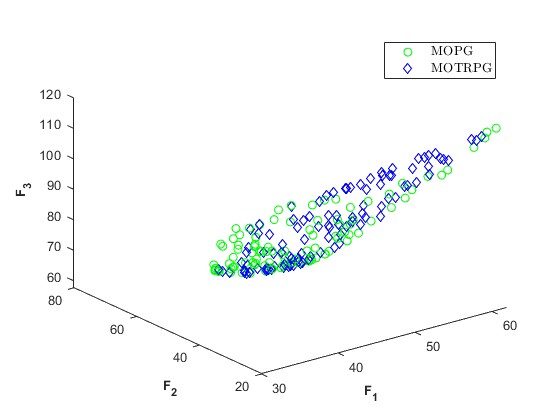}
        \caption{Approximate Pareto fronts of (MOLS)}
        \label{fig2}
\end{figure}

Multi-objective least square problems can be solved using Algorithm \ref{trpxm1} using 
\bea
f_j(x)&=&\frac{1}{2}\|A_jx-b_j\|^2 \lb{ls1a}\\
g_j(x)&=&\frac{\nu_j}{2} \|x\|_1. \lb{ls1b}
\eea
We have constructed a three objective optimization problem $$\underset{x\in\mathbb{R}^3}{\min}~~\left(\frac{1}{2} \|A_1x-b_1\|^2+\frac{\nu_1}{2} \|x\|_1, \frac{1}{2} \|A_2x-b_2\|^2+\frac{\nu_2}{2} \|x\|_1, \frac{1}{2} \|A_2x-b_2\|^2+\frac{\nu_3}{2} \|x\|_1 \right) $$
where entries of $(A_j)\in \mathbb{R}^{10\times 3}$ and $b_j\in\mathbb{R}^{10}$ are generated using uniform distribution in $[0,5]$ and $[0,10]$ respectively. For nonsmooth functions, we have considered $\nu_1=0.30$, $\nu_2=1.06$, $\nu_3=1.84$. A set of 100 initial approximation uniformly distributed in $(-1,-1,-1)^T$ and $(1,1,1)^T$ is used to find approximate Pareto front. Approximate Pareto front generated by $(MOTRPG)$ and $(MOPG)$ if provided in Figure \ref{fig2}.
$ $\\

MOTRPG is compared with MONPG and MOPG. Each algorithm is executed on a set of test problems. Smooth functions ($f_j$) of the test problems are collected from different sources and non smooth functions ($g_j$) are either collected from \cite{mat6} ($gA$, $gB$, $gC$, and $gF$ of Section 5 in \cite{mat6}) or (\ref{ls1b}) with $\nu\in\mathbb{R}^m$ uniformly generated in $0^m$ and $2^m$. Details of test problems are provided in Table \ref{table1}. In this table $m$, $n$ denotes the number of objective functions and dimension of $x$ respectively. A set of 100 initial points uniformly distributed between $lb$ and $ub$ is considered to generate approximate Pareto front.
\begin{table}[!http]
\begin{center}
\scriptsize
\begin{tabular}
{|c|c|c|c|c|c|}\hline
Sl. No& $(m, n)$& $f$ & $g$ &$lb^T$ &$ub^T$\\ \hline
\rownumber&(2,2)& Ansary1 (\cite{mat6})& $gA$  &$(-3,-3)$& $(7,7)$ \\ \hline
\rownumber&(2,2)& Ansary1 (\cite{mat6})& $gB$  &$(-3,-3)$& $(7,7)$\\ \hline
\rownumber&(2,2)&AP2 (\cite{mat1})& gA &$(-5,-5)$&$(5,5)$\\ \hline
\rownumber &(2,2)&AP2 (\cite{mat1})& gB &$(-5,-5)$&$(5,5)$\\ \hline
\rownumber&(2,2)&BK1 (\cite{hub1})& gA &$(-5,-5)$&$(7.5,7.5)$\\ \hline
\rownumber &(2,2)&BK1 (\cite{hub1})& gB &$(-5,-5)$&$(7.5,7.5)$\\ \hline
\rownumber &(3,3)&FDS (\cite{flg1})& gC &$(-2,-2,-2)$&$(4,4,4)$\\ \hline
\rownumber &(3,5)&FDS (\cite{flg1})& (\ref{ls1b})&$(-2,\dots,-2)$&$(2,\dots,2)$\\ \hline
\rownumber &(3,8)&FDS (\cite{flg1})& (\ref{ls1b})&$(-2,\dots,-2)$&$(2,\dots,2)$\\ \hline
\rownumber &(3,2)&IKK1 (\cite{hub1})& (\ref{ls1b})&$(-2,-2)$&$(3,3)$\\ \hline
\rownumber &(2,2)&Jin1 (\cite{jin1})& gA&$(-3,-3)$&$(5,5)$\\ \hline
\rownumber &(2,2)&Jin1 (\cite{jin1})& gB &$(-3,-3)$&$(5,5)$\\ \hline
\rownumber &(2,4)&Jin1 (\cite{jin1})& gF&$(-5,\dots,-5)$&$(10,\dots,10)$\\ \hline
\rownumber &(2,10)&Jin1 (\cite{jin1})& (\ref{ls1b}) &$(-5,\dots,-5)^T$&$(5,\dots,5)^T$\\ \hline
\rownumber  &(2,2)&Lovison1 (\cite{lovison1})& gA &$(-3,-3)$&$(5,5)$\\ \hline
\rownumber &(2,2)&Lovison1 (\cite{lovison1})& gB&$(-3,-3)$&$(5,5)$\\ \hline
\rownumber &(2,2)&Lovison4 (\cite{lovison1})& gB&$(-10,-5)$&$(5,5)$\\ \hline
\rownumber &(2,2)&LRS1 (\cite{hub1})& gB&$(-50,-50)$&$(50,50)$\\ \hline
\rownumber &(2,2)&LRS1 (\cite{hub1})& gA&$(-50,-50)$&$(50,50)$\\ \hline
\rownumber &(3,3)&MOLS1 (\ref{ls1a})& (\ref{ls1b}) &$(-1,-1,-1)$&$(1,1,1)$\\ \hline
\rownumber &(3,1)&MHHM1 (\cite{hub1})& (\ref{ls1b}) &$-4$&$4$\\ \hline
\rownumber &(3,2)&MHHM2 (\cite{hub1})& (\ref{ls1b}) &$(-4,-4)$&$(4,4)$ \\ \hline 
\rownumber &(2,1)&MOP1 (\cite{hub1})& (\ref{ls1b})&$-100$&$100$\\ \hline 
\rownumber &(3,2)&MOP7 (\cite{hub1})& (\ref{ls1b})&$(-4,-4)$&$(4,4)$\\ \hline 
\rownumber &(2,2)&MS1(\cite{martin1})& gA&$(-2,-2)$&$(2,2)$\\ \hline
\rownumber &(2,2)&MS1 (\cite{martin1})& gB&$(-2,-2)$&$(2,2)$\\ \hline
\rownumber &(2,4)&MS2(\cite{martin1})& gG&$(-2,\dots,-2)$&$(2,\dots,2)$\\ \hline
\rownumber &(2,10)&MS2 (\cite{martin1})&  (\ref{ls1b})&$(-2,\dots,-2)$&$(2,\dots,2)$\\ \hline
\rownumber &(3,3)&SDD1 (\cite{sdd1})& gC&$(-2,-2,-2)$&$(2,2,2)$\\ \hline
\rownumber &(3,10)&SDD1 (\cite{sdd1}) &  (\ref{ls1b})&$(-2,\dots,-2)$&$(2,\dots,2)$\\ \hline
\rownumber&(2,2)&SP1 (\cite{hub1})& gA&$(-1,-1)$&$(5,5)$\\ \hline
\rownumber &(2,2)&SP1 (\cite{hub1})& gB&$(-1,-1)$&$(5,5)$\\ \hline 
\rownumber&(2,2)&SSFY1 (\cite{hub1})& gA&$(-50,-50)$&$(50,50)$\\ \hline
\rownumber &(2,2)&SSFY1 (\cite{hub1})& gB&$(-50,-50)$&$(50,50)$\\ \hline 
\rownumber &(3,2)&VFM1 (\cite{hub1})&  (\ref{ls1b}) &$(-2,-2)$&$(2,2)$ \\ \hline
\rownumber&(2,2)&VU1 (\cite{hub1})& gA&$(-3,-3)$&$(3,3)$\\ \hline
\rownumber &(2,2)&VU1 (\cite{hub1})& gB&$(-3,-3)$&$(3,3)$\\ \hline
\rownumber&(2,2)&VU2 (\cite{hub1})& gA&$(-3,-3)$&$(3,3)$\\ \hline
\rownumber &(2,2)&VU2(\cite{hub1})& gB&$(-3,-3)$&$(3,3)$\\ \hline 
\rownumber &(3,3)&ZLT1 (\cite{hub1})& gC &$(-100,-100,-100)$&$(100,100,100)$\\ \hline 
\rownumber &(3,10)&ZLT1 (\cite{hub1})&  (\ref{ls1b}) &$(-5,\dots,-5)$&$(5,\dots,5)$\\ \hline
\end{tabular}
\caption{Details of test problems}
\lb{table1}
\end{center}
\end{table}
$ $\\

Performance profiles using `purity metric', `$\Delta$ and $\Gamma$ spread metrics', `hypervolome metric', and 'total function evaluations' are used to compare different methods. Detailed explanation of performance profile as well as `purity metric', `$\Delta$ and $\Gamma$ spread metrics', and `hypervolome metric' are provided in \cite{mat6,mat2,mat4,mat5,flg3}. Total function evaluation of an algorithm for a test problem is based on total function, gradient, and Hessian (if required) evaluations. Suppose $\# f$, $\# \nabla f $, and $\# \nabla^2 f$ denote the number of function, gradient, and Hessian evaluations required to solve a test problem. We have used forward difference formula to compute approximate $\nabla f(x)$ and $\nabla^2 f(x)$. For MOTRPG and MOPG, total function evaluations is $$\# ~fun= \# f+ n *\# \nabla f.$$ For MONPG, total function evaluations is $$\# ~fun= \# f+ n *\# \nabla f+\frac{n*(n+1)}{2} * \nabla^2 f.$$
$ $\\

\begin{figure}[!htbp]
    \centering
    \begin{subfigure}[b]{.47\textwidth}
        \centering
        \includegraphics[height=2.5cm,width=\linewidth]{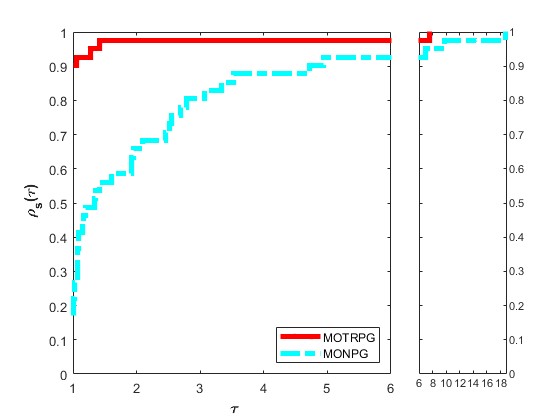}
        \caption{Between MOTRPG and MONPG}
        \label{pua}
    \end{subfigure}
    \begin{subfigure}[b]{.47\textwidth}
        \centering
        \includegraphics[height=2.5cm,width=\linewidth]{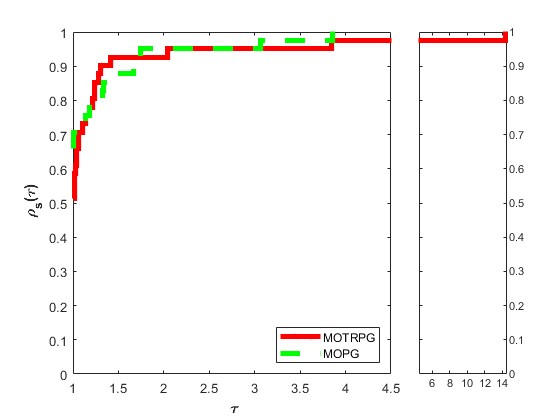}
        \caption{Between MOTRPG and MOPG}
        \label{pub}
    \end{subfigure}
\label{fig_pu}
\caption{Performance profiles using purity metric}
\end{figure}

\begin{figure}[!htbp]
    \centering
    \begin{subfigure}[b]{.47\textwidth}
        \centering
        \includegraphics[height=2.5cm,width=\linewidth]{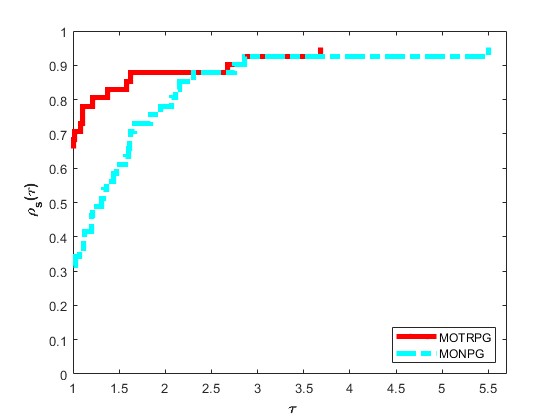}
        \caption{Between MOTRPG and MONPG}
        \label{gamma}
    \end{subfigure}
    \begin{subfigure}[b]{.47\textwidth}
        \centering
        \includegraphics[height=2.5cm,width=\linewidth]{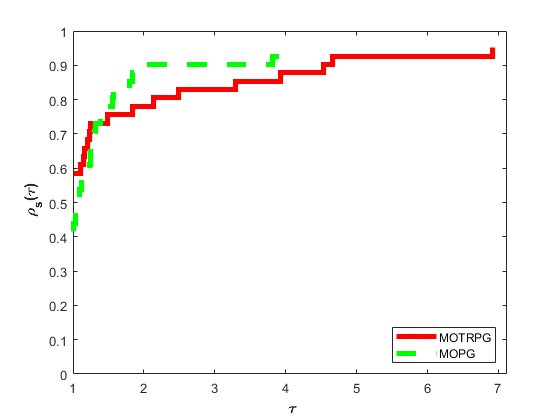}
        \caption{Between MOTRPG and MOPG}
        \label{gammb}
    \end{subfigure}
\caption{Performance profiles using $\Gamma$-spread metric}
\end{figure}

\begin{figure}[!htbp]
    \centering
    \begin{subfigure}[b]{.47\textwidth}
        \centering
        \includegraphics[height=2.5cm,width=\linewidth]{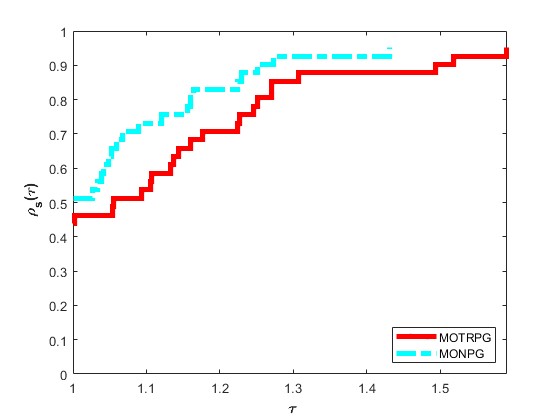}
        \caption{Between MOTRPG and MONPG}
        \label{deltaa}
    \end{subfigure}
    \begin{subfigure}[b]{.47\textwidth}
        \centering
        \includegraphics[height=2.5cm,width=\linewidth]{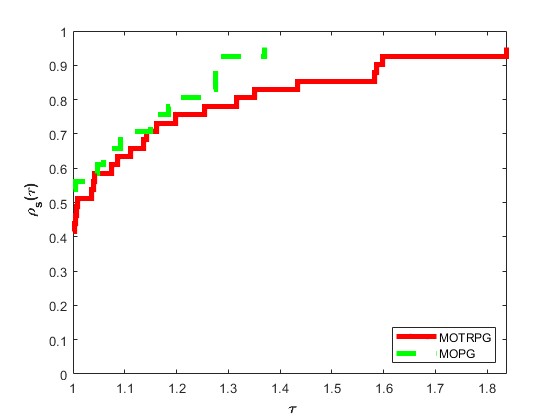}
        \caption{Between MOTRPG and MOPG}
        \label{deltab}
    \end{subfigure}
\caption{Performance profiles using $\Delta$-spread metric}
\end{figure}

\begin{figure}[!htbp]
    \centering
    \begin{subfigure}[b]{.47\textwidth}
        \centering
        \includegraphics[height=2.5cm,width=\linewidth]{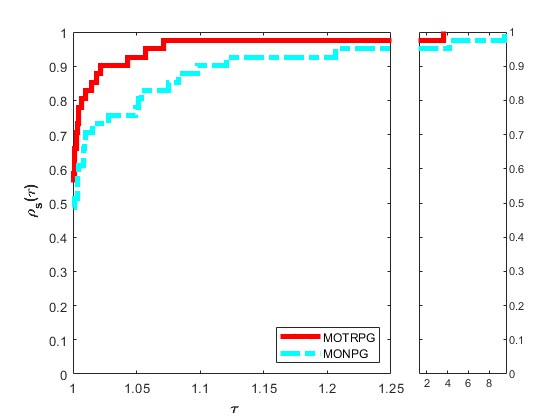}
        \caption{Between MOTRPG and MONPG}
        \label{hva}
    \end{subfigure}
    \begin{subfigure}[b]{.47\textwidth}
        \centering
        \includegraphics[height=2.5cm,width=\linewidth]{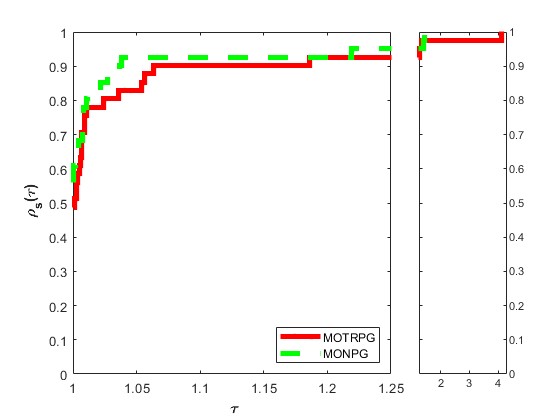}
        \caption{Between MOTRPG and MOPG}
        \label{hvb}
    \end{subfigure}
\caption{Performance profiles using hypervolume metric}
\end{figure}

\begin{figure}[!htbp]
    \centering
    \begin{subfigure}[b]{.47\textwidth}
        \centering
        \includegraphics[height=2.5cm,width=\linewidth]{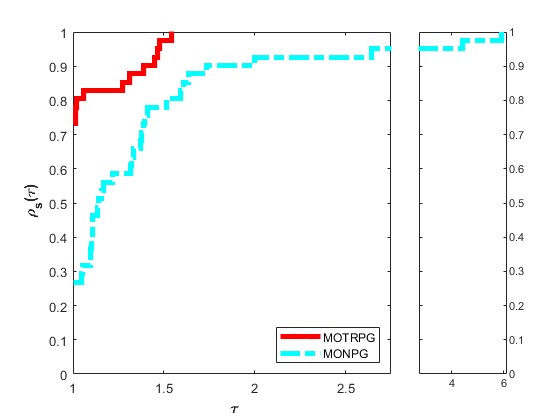}
        \caption{Between MOTRPG and MONPG}
        \label{fca}
    \end{subfigure}
    \begin{subfigure}[b]{.47\textwidth}
        \centering
        \includegraphics[height=2.5cm,width=\linewidth]{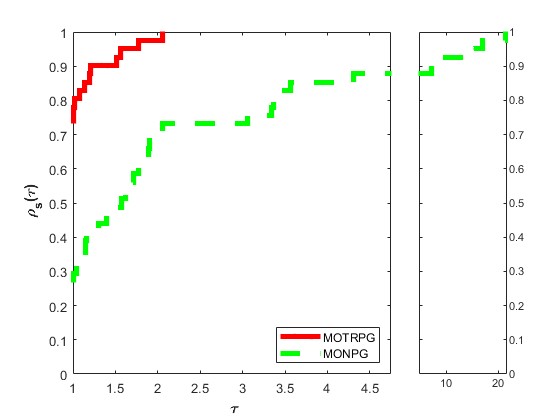}
        \caption{Between MOTRPG and MOPG}
        \label{fcb}
    \end{subfigure}
    \label{fig_fc}
\caption{Performance profiles using \# fun}
\end{figure}

Performance profile between MOTRPG and MONPG  using purity matrix is provided in Figure \ref{pua}.  Performance profile between MOTRPG and MOPG  using purity matrix is provided in Figure \ref{pub}. Figure \ref{gamma} and \ref{gammb} represent the performance profiles using $\Gamma$-spread metric. Similarly performance profiles using $\Delta$-spread metric are provided in Figure \ref{deltaa} and \ref{deltab}.  Figure \ref{hva} and \ref{hvb} represent the performance profiles using  hypervolume metric. Last bu not least, performance profile unsing function evaluations are provided in Figures \ref{fca} and \ref{fcb}.
$ $\\
From performance profile figures (Figures \ref{pua}-\ref{fcb}), one can observe that MONPG performs better than MONPG and MOPG in purity, $\Gamma$-spread metric, and total function evaluations in most cases. Other metrics, MOTRPG performs average with respect to other methods. 
\section{Conclusion}
In this paper, we have developed a trust region proximal gradient method nonlinear convex multi-objective optimization problems. This method is free from any kind of priori chosen parameters or ordering information of objective function. Global convergence of the proposed method is justified under some mild assumptions. However, the proposed method is restricted to convex multi-objective optimization problems only. In future we want to extend the results to nonconvex multi-objective optimization problems. We have used multi-start technique to generate approximate Pareto fronts, which fails to generate a well distributed approximate Pareto front in some cases. Some initial point selection technique can be developed to generate a well distributed approximate Pareto front, which is left as a future scope of this method. 
\section*{Declarations}
{\bf Ethics approval and consent to participate:} Not applicable.\\
{\bf Consent for publication:} The author declares no consent for publication.\\
{\bf Data availability statement:} This paper does not involve any associated data.\\
{\bf Competing interest:} The author declares no conflict interest.\\
{\bf Funding information:} No\\
{\bf Authors' contributions:} Single author has contributed everything.\\
{\bf Acknowledgements:} The author is thankful to IIT Jodhpur for giving research facilities.

\end{document}